\documentclass[11pt]{amsart}
\usepackage[parfill]{parskip} 
\usepackage{amsmath, amsxtra, amsthm, amssymb,dsfont,eucal,stmaryrd}
\usepackage[all]{xy}
\usepackage{color}
\usepackage[colorlinks=true, pdfstartview=FitV, linkcolor=blue, citecolor=blue, urlcolor=blue]{hyperref}

\usepackage{graphicx}

\numberwithin{equation}{section}
\swapnumbers

\newcommand{\nc}{\newcommand}
\nc{\rc}{\renewcommand}

\rc{\b}{\mathbb}
\rc{\c}{\mathcal}
\nc{\on}{\operatorname}
\nc{\tn}{\textnormal}

\nc{\bA}{\b A}
\nc{\bB}{\b B}
\nc{\bC}{\b C}
\nc{\bD}{\b D}
\nc{\bE}{\b E}
\nc{\bF}{\b F}
\nc{\bG}{\b G}
\nc{\bH}{\b H}
\nc{\bI}{\b I}
\nc{\bJ}{\b J}
\nc{\bK}{\b K}
\nc{\bL}{\b L}
\nc{\bM}{\b M}
\nc{\bN}{\b N}
\nc{\bO}{\b O}
\nc{\bP}{\b P}
\nc{\bQ}{\b Q}
\nc{\bR}{\b R}
\nc{\bS}{\b S}
\nc{\bT}{\b T}
\nc{\bU}{\b U}
\nc{\bV}{\b V}
\nc{\bW}{\b W}
\nc{\bX}{\b X}
\nc{\bY}{\b Y}
\nc{\bZ}{\b Z}

\nc{\cA}{\c A}
\nc{\cB}{\c B}
\nc{\cC}{\c C}
\nc{\cD}{\c D}
\nc{\cE}{\c E}
\nc{\cF}{\c F}
\nc{\cG}{\c G}
\nc{\cH}{\c H}
\nc{\cI}{\c I}
\nc{\cJ}{\c J}
\nc{\cK}{\c K}
\nc{\cL}{\c L}
\nc{\cM}{\c M}
\nc{\cN}{\c N}
\nc{\cO}{\c O}
\nc{\cP}{\c P}
\nc{\cQ}{\c Q}
\nc{\cR}{\c R}
\nc{\cS}{\c S}
\nc{\cT}{\c T}
\nc{\cU}{\c U}
\nc{\cV}{\c V}
\nc{\cW}{\c W}
\nc{\cX}{\c X}
\nc{\cY}{\c Y}
\nc{\cZ}{\c Z}

\nc{\fA}{{\mathfrak A}}
\nc{\fB}{{\mathfrak B}}
\nc{\fC}{{\mathfrak C}}
\nc{\fD}{{\mathfrak D}}
\nc{\fE}{{\mathfrak E}}
\nc{\fF}{{\mathfrak F}}
\nc{\fG}{{\mathfrak G}}
\nc{\fH}{{\mathfrak H}}
\nc{\fI}{{\mathfrak I}}
\nc{\fJ}{{\mathfrak J}}
\nc{\fK}{{\mathfrak K}}
\nc{\fL}{{\mathfrak L}}
\nc{\fM}{{\mathfrak M}}
\nc{\fN}{{\mathfrak N}}
\nc{\fO}{{\mathfrak O}}
\nc{\fP}{{\mathfrak P}}
\nc{\fQ}{{\mathfrak Q}}
\nc{\fR}{{\mathfrak R}}
\nc{\fS}{{\mathfrak S}}
\nc{\fT}{{\mathfrak T}}
\nc{\fU}{{\mathfrak U}}
\nc{\fV}{{\mathfrak V}}
\nc{\fW}{{\mathfrak W}}
\nc{\fZ}{{\mathfrak Z}}
\nc{\fX}{{\mathfrak X}}
\nc{\fY}{{\mathfrak Y}}
\nc{\fa}{{\mathfrak a}}
\nc{\fb}{{\mathfrak b}}
\nc{\fc}{{\mathfrak c}}
\nc{\fd}{{\mathfrak d}}
\nc{\fe}{{\mathfrak e}}
\nc{\ff}{{\mathfrak f}}
\nc{\fg}{{\mathfrak g}}
\nc{\fh}{{\mathfrak h}}
\nc{\fiI}{{\mathfrak i}}  
\nc{\ffi}{{\mathfrak i}}  
\nc{\fj}{{\mathfrak j}}
\nc{\fk}{{\mathfrak k}}
\nc{\fl}{{\mathfrak{l}}}
\nc{\fm}{{\mathfrak m}}
\nc{\fn}{{\mathfrak n}}
\nc{\fo}{{\mathfrak o}}
\nc{\fp}{{\mathfrak p}}
\nc{\fq}{{\mathfrak q}}
\nc{\fr}{{\mathfrak r}}
\nc{\fs}{{\mathfrak s}}
\nc{\ft}{{\mathfrak t}}
\nc{\fu}{{\mathfrak u}}
\nc{\fv}{{\mathfrak v}}
\nc{\fw}{{\mathfrak w}}
\nc{\fz}{{\mathfrak z}}
\nc{\fx}{{\mathfrak x}}
\nc{\fy}{{\mathfrak y}}

\nc{\al}{{\alpha }}
\nc{\be}{{\beta }}
\nc{\ga}{{\gamma }}
\nc{\de}{{\delta }}
\nc{\del}{{\partial }}
\nc{\vep}{{\varepsilon }}
\nc{\ep}{{\epsilon }}

\nc{\ze}{{\zeta }}
\nc{\et}{{\eta }}
\rc{\th}{{\theta }}

\nc{\vth}{{\vartheta }}

\nc{\io}{{\iota }}
\nc{\ka}{{\kappa }}
\nc{\la}{{\lambda }}
\nc{\vrho}{{\varrho}}
\nc{\si}{{\sigma }}
\nc{\ups}{{\upsilon }}
\nc{\vphi}{{\varphi }}
\nc{\om}{{\omega }}

\nc{\Ga}{{\Gamma }}
\nc{\De}{{\Delta }}
\nc{\nab}{{\nabla}}
\nc{\Th}{{\Theta }}
\nc{\La}{{\Lambda }}
\nc{\Si}{{\Sigma }}
\nc{\Ups}{{\Upsilon }}
\nc{\Om}{{\Omega }}

\nc{\Spec}{\on{Spec}}
\nc{\id}{\on{id}}
\nc{\inv}{ ^{-1}}
\nc{\su}{\subset}
\nc{\ot}{\otimes}
\nc{\un}{\underline}
\nc{\ov}{\overline}
\nc{\uu}{\mathds{1}}

\nc{\shom}{\cH om}
\nc{\Gr}{\cG\tn{r}}

\nc{\se}{\section}
\nc{\sse}{\subsection}
\nc{\ssse}{\subsubsection}

\nc{\wt}{\widetilde}
\nc{\slc}{\fs \fl_2\bC}
\nc{\gr}{\tn{gr}}
\nc{\lan}{\langle}
\nc{\ran}{\rangle}
\nc{\ord}{\tn{ord}}
\nc{\wh}{\widehat}
\nc{\sst}{\scriptstyle}
\nc{\sss}{\scriptscriptstyle}

\nc{\qchoose}[2]{\left[\!\! \begin{array}{c} #1 \\ #2 \end{array} \!\!\right]_q}
\nc{\bn}{\bf n}
\nc{\emp}{\emptyset}
\nc{\plat}{\tn{plat}}
\nc{\spr}{\tn{spr}}
\nc{\bk}{\backslash}
\nc{\md}[1]{\tn{ (mod }#1)}

\nc{\lf}{\lfloor}
\nc{\rf}{\rfloor}
\nc{\llf}{\left\lfloor}
\nc{\rrf}{\right\rfloor}

\title{A combinatorial proof of strict unimodality for $q$-binomial coefficients}

\author{Vivek Dhand}

\begin{document}

\maketitle
\thispagestyle{empty}

\begin{abstract}
I.\,Pak and G.\,Panova recently proved that the $q$-binomial coefficient ${m+n \choose m}_q$ is a strictly unimodal polynomial in $q$ for $m,n \geq 8$, via the representation theory of the symmetric group.  We give a direct combinatorial proof of their result by characterizing when a product of chains is strictly unimodal and then applying O'Hara's structure theorem for the partition lattice $L(m,n)$.  In fact, we prove a stronger result: if $m, n \geq 8d$, and $2d \leq r \leq mn/2$, then the $r$-th rank of $L(m,n)$ has at least $d$ more elements that the next lower rank.
\end{abstract}

\se{Introduction}

Recall that the lattice $L(m,n)$ consists of integer partitions whose Young diagrams fit inside an $(m \times n)$-rectangle, ordered by inclusion:
\[ L(m,n) = \{ (\la_1, \dots, \la_m) \mid n \geq \la_1 \geq \dots \la_m \geq 0  \}. \]
The rank generating function for $L(m,n)$ is the $q$-binomial coefficient:
\[ {m+n \choose m}_q = \prod_{i = 1}^n \frac{1-q^{m+i}}{1-q^i} = \sum_{r = 0}^{mn} p_r(m,n) q^r. \]
It is easy see that complementary ranks have the same size: $p_r(m,n) = p_{mn-r}(m,n)$.
Sylvester was the first to prove the {\em unimodality} property \cite{Syl}:
\[ p_0(m,n) \leq \dots \leq p_{\lf mn/2 \rf}(m,n). \]
Several other proofs have been discovered over the years using many different techniques, e.g.~\cite{P,S,Ze}.  In particular, O'Hara gave a purely combinatorial proof by decomposing the underlying ranked set of $L(m,n)$ into a centered disjoint union of products of chains \cite{O}.  We use a slightly modified version of O'Hara's theorem to prove the following:

\sse{Theorem} If $m, n \geq 8d$ and $2d \leq r \leq mn/2$, then:
\[ p_r(m,n) - p_{r-1}(m,n) \geq d. \]

If $d=1$, we recover the strict unimodality theorem of Pak-Panova \cite{PP}:

\sse{Theorem}  If $m,n \geq 8$ and $2 \leq r \leq mn/2$, then:
\[ p_r(m,n) - p_{r-1}(m,n) > 0. \] 




 In what follows we will ignore the partial order on $L(m,n)$ and simply work with the underlying ranked set.   The following is a restatement of O'Hara's structure theorem \cite{D2,O}:

\sse{Theorem} There is a decomposition of $L(m,n)$ into centered rank-symmetric unimodal subsets $Q_m(d_0, \dots, d_k)$, where:
\[ k = \lf m/2 \rf \quad \tn{ and } \quad n = d_0 + 2d_1 + \dots + (k+1)d_k, \]
and there exist isomorphisms of ranked sets:
\[  Q_m(d_0, \dots, d_k) \simeq L(r,\ell_m(d_0, \dots, d_k)) \times Q_{m-2r}(d_r, \dots, d_k) \]
where:
\[ r = 1 + \min\{ j \mid d_j > 0\}  \quad \tn{ and } \quad \ell_m(d_0, \dots, d_k) = \sum_{j = 0}^k (m-2j) d_j. \]

\sse{Remark} There are several equivalent ways to define the sets $Q_m(d_0, \dots, d_k)$.  For example, they are the level sets of certain tropical polynomials \cite{D1}:
\[ f_{m,r}(\la) =  \min_{0\leq t_0 \leq \dots \leq t_{m-2r} \leq r}\sum_{j = 0}^{m-2r} a_{2 t_j + j}, \]
where $\la_0 = n, \la_{m+1} = 0$, and $a_i = \la_i - \la_{i+1}$.  We have the formula:
\[ Q_m(d_0, \dots, d_k) = \{ \la \in L(m,n) \mid f_{m,r}(\la) = \sum_{j = r}^k (j + 1 - r)d_j \tn{ for } 1 \leq r \leq k \}. \]
In particular, it follows that $Q_m(d_0, \dots, d_k)$ has a unique minimal element, and its rank is equal to:
\[ \sum_{j = 0}^k j(j+1) d_j. \]

\sse{Remark} Our proof involves choosing subsets $Q_m(d_0, \dots, d_k) \su L(m,n)$ whose strictly unimodal ranges provide a covering of the desired interval.  F. Zanello has given a similar combinatorial proof \cite{Za}.

\se{Strict unimodality for chain products}

Let $P$ be a ranked poset $P$ of length $n$.  Let $p_i$ denote the size of the $i$-th rank of $P$.  We say that $P$ is {\em rank-symmetric} if $p_i = p_{n-i}$ for all $0 \leq i \leq n$.  We say that a rank-symmetric poset $P$ is {\em strictly unimodal} if:
\[ p_0 < \dots < p_{\lf n/2\rf}. \]
In this section, we determine when a product of chains is strictly unimodal. We will apply this result in the next section to prove the theorem. 

Given a non-negative integer $a$, let $[a] = \{0 < \dots < a\}$ denote a chain of length $a$.  

\sse{Lemma}  Let $P = [a_1] \times \dots \times [a_n]$, where $a_1 \geq \dots \geq a_n \geq 0$.  Let $C$ be a shortest chain in a symmetric chain decomposition of $P$.

(1) $P$ is strictly unimodal from rank $0$ up to the lowest rank of $C$.

(2) The length of $C$ is equal to:
\[ \ell(C) =  \max(a_1 - (a_2 + \dots + a_n),\ep) \]
where $\ep = 0$ (resp.~$\ep = 1$) if $a_1 + \dots + a_n$ is even (resp.~odd).

(3) $P$ is strictly unimodal if and only if:
\[ a_1 \leq a_2 + \dots + a_n + 1. \]

\begin{proof}
We will simultaneously prove these statements by induction on $n$.  Consider the base case $n = 2$.  There is a well-known symmetric chain decomposition:
\[ [a_1] \times [a_2] \simeq \bigsqcup_{i = 0}^{a_2} [a_1 + a_2 - 2i] \]
where the lowest rank of the $i$-th chain is $i$.  Since there is a new symmetric chain starting at each rank from $0$ to $a_2$, we see that $[a_1] \times [a_2]$ is strictly unimodal from rank $0$ up to $a_2$, which is equal to the lowest rank of the shortest chain.  The length of the shortest chain is:
\[ a_1 + a_2 - 2a_2 = a_1 - a_2, \]
so $[a_1] \times [a_2]$ is strictly unimodal if and only if $a_1 - a_2 \leq 1$.

For $n \geq 3$, let $R = [a_2] \times \dots \times [a_n]$ so:
\[ P = [a_1] \times R. \]
Let $D$ be a shortest chain in a symmetric chain decomposition of $R$, and let $r$ denote the lowest rank of $D$.  By induction, $R$ is strictly unimodal from rank $0$ up $r$ and the length of $D$ is equal to:
\[ \ell(D) = \max(a_2 - (a_3 + \dots + a_n),\ep') \]
where $\ep' = 0$ (resp.~$\ep' = 1$) if $a_2 + \dots + a_n$ is even (resp.~odd).  In any case, we have:
\[ a_1 \geq a_2 \geq \ell(D).  \]
Now, the length of the {\em longest} chain in $R$ is $a_2 + \dots + a_n$, and the lengths of all the symmetric chains in $R$ have the same parity.

If $a_1 \leq a_2 + \dots + a_n$, then the smallest difference between $a_1$ and the lengths of all the symmetric chains in $R$ is equal to $0$ or $1$, depending on the parity of $a_1 + \dots + a_n$.  In other words,  the length of a shortest chain $C$ in any SCD of $P$ will be $0$ or $1$.  On the other hand, if $a_1 > a_2 + \dots + a_n$, then:
\[ \ell(C) = a_1 - (a_2 + \dots + a_n) \]
by our formula from the base case.  Therefore, we obtain:
\[ \ell(C) =  \max(a_1 - (a_2 + \dots + a_n),\ep). \]
Note that $P$ is strictly unimodal from rank $0$ up to the lowest rank of $C$ because $P$ is a disjoint union over products of the form $[a_1] \times D'$, where $D'$ runs over the symmetric chains in a fixed SCD of $R$.  Finally, $P$ is strictly unimodal if and only if $\ell(C) \leq 1$, and the result follows.
\end{proof}

\se{Proof of the theorem}

\sse{Theorem} If $m, n \geq 8d$ and $2d \leq r \leq mn/2$, then:
\[ p_r(m,n) - p_{r-1}(m,n) \geq d. \]

\begin{proof} {\bf Step 1:} Consider $Q_m(d_0, d_1, 0, \dots, 0) \su L(m,n)$ where $d_0, d_1 > 0$.  This ranked set is isomorphic to a product of two chains:
\[ [md_0 + (m-2)d_1] \times [(m-2)d_1]. \]
The lowest rank element of this set has rank $2d_1$.  Note that $n = d_0 + 2d_1$ so $n > 2d_1$. We know that a product of two chains of the form $[a+b] \times [b]$ is strictly unimodal from rank $0$ to rank $b$, so $Q_m(d_0, d_1, 0, \dots, 0)$ is strictly unimodal from rank $2d_1$ up to rank:
\[  2d_1 + (m-2)d_1 = md_1. \]
Therefore, we see that:
\[ p_r(m,n) - p_{r-1}(m,n) \geq d \]
for all $2d \leq r \leq mn/2 - dm$ and $n > 2d$.

{\bf Step 2:} Consider $Q_m(d_0,d_1,d_2, 0, \dots, 0) \su L(m,n)$ where $d_0, d_1, d_2 > 0$.  This ranked set is isomorphic to a product of three chains:
\[  [md_0 + (m-2)d_1 + (m-4)d_2] \times [(m-2)d_1 + (m-4)d_2] \times [(m-4)d_2]. \]
The lowest rank element of this set has rank $2d_1 + 6d_2$.  We know that the product of three chains of the form $[a+b+c] \times [b+c] \times [c]$ is strictly unimodal if and only if $a \leq c + 1$, so $Q_m(d_0, d_1, d_2, 0, \dots, 0)$ is strictly unimodal if and only if:
\[ md_0 \leq (m-4)d_2 + 1. \]
Let $N_3(m,n)$ denote the number of strictly unimodal subsets in $L(m,n)$ of the above type:
\[ N_3(m,n) = \{ (d_0, d_1, d_2) \in \bZ_{>0}^3 \mid d_0 + 2d_1 + 3d_2 = n, \; md_0 \leq (m-4)d_2 + 1\}. \] 
Now let us calculate the largest possible value of $d_2$, depending on $n \md{3}$:
\[ n = 0 \md{3} \implies (d_0, d_1, d_2 ) = (1,1, (n-3)/3). \]
\[ n = 1 \md{3} \implies (d_0, d_1, d_2 ) = (2,1, (n-4)/3). \]
\[ n = 2 \md{3} \implies (d_0, d_1, d_2 ) = (1,2, (n-5)/3). \]
Therefore, the largest possible value of $2d_1 + 6d_2$ is equal to $2n - 4$ if $n = 0 \md{3}$ and $2n - 6$ otherwise.  

It follows that:
\[ p_r(m,n) - p_{r-1}(m,n) \geq N_3(m,n) \]
for all $2n-4 \leq r \leq mn/2$.

{\bf Step 3:}  We see that we are reduced to satisfying the following inequalities:
\[ mn/2- dm \geq 2n-4 \quad \tn{and} \quad N_3(m,n) \geq d. \]
The first inequality is equivalent to:
\[ m \geq \frac{4n-8}{n-2d}. \]
If $d \geq 1$ and $m, n \geq 8d$, then the above inequality does hold:
\[ m(n-2d) \geq 8d(n-2d) = d(8n - 16d) \geq 8n-16d \geq 4n + 32d - 16d = 4n + 16d \geq 4n-8. \]
Let us now estimate a lower bound for $N_3(m,n)$.  For a fixed value of $d_0$, we find that:
\[ 3d_2 = n - d_0 - 2d_1 \leq n - 2 - d_0. \]
Given a possible solution $(d_1, d_2)$, note that the next solution is $(d_1 + 3, d_2 - 2)$, so all the possible values of $d_2$ must have the same parity mod $2$.
Therefore:
\[ \frac{m d_0 -1}{m-4} \leq d_2 \leq \frac{n - 2 -d_0}{3} \]
where all the values of $d_2$ must have the same parity mod $2$.  So the number of solutions is at least:
\[ \frac{1}{2} \left( \frac{n - 2 - d_0}{3} - \frac{md_0 -1}{m-4} \right). \]
Summing over the allowed values of $d_0$, we obtain:
\[ N_3(m,n) \geq \frac{1}{2} \sum_{d_0} \left( \frac{n - 2 - d_0}{3} - \frac{md_0 -1}{m-4} \right). \]
Note that, if $m \geq 8$, then:
\[ md_0 -1 < 2d_0(m-4) = 2md_0 - 8d_0 \]
because:
\[ md_0 \geq 8d_0 > 8 d_0 -1. \]
Therefore:
\[  \frac{1}{2}\left( \frac{n - 2 - d_0}{3} - \frac{md_0 -1}{m-4} \right) > \frac{1}{2}\left( \frac{n - 2 - d_0}{3} - 2d_0 \right) = \frac{n - 2 - 7d_0}{6}. \]
If $d \geq 3$, then $m,n \geq 24$ and by considering the first three terms we obtain:
\[ N_3(m,n) \geq \frac{n - 9}{6} + \frac{n - 16}{6} + \frac{n-23}{6} = \frac{n-16}{2} \geq 4d-8 \geq d \]
for all $d \geq 3$.  It remains to prove the theorem for $1 \leq d \leq 2$. 

{\bf Step 4:}  Consider $Q_m(d_0,d_1,d_2, d_3, 0, \dots, 0) \su L(m,n)$ where $d_0, d_1, d_2,d_3 > 0$.  This ranked set is isomorphic to a product of four chains:
\[ [\ell_0] \times  [\ell_1] \times  [\ell_2] \times  [\ell_3], \]
where $\ell_3 = (m-6)d_3$, $\ell_2 = \ell_3 + (m-4)d_2$, $\ell_1 = \ell_2 + (m-2)d_1$, and $\ell_0 = \ell_1 + md_0$.  It is strictly unimodal if and only if:
\[ \ell_0 - \ell_1 \leq \ell_2 + \ell_3  + 1, \]
which translates to the condition that:
\[ md_0 \leq (m-4)d_2 + 2(m-6)d_3. \]
If $d_2, d_3 \geq d_0$ we get:
\[ (m-4)d_2 + 2(m-6)d_3 \geq d_0 (3m - 16) \geq md_0  \]
for all $m \geq 8$.  In other words, $Q_m(1,d_1,d_2, d_3, 0, \dots, 0) \su L(m,n)$ is strictly unimodal if there exist $x,y,z \geq 0$ such that $2x+3y+4z = n-10$.  The standard generating function argument shows that the existence of at least one solution for all $n \geq 12$ and at least two solutions for all $n \geq 16$.  So the $d = 2$ case is finished, and for $d = 1$ we can check the remaining cases $8 \leq m,n \leq 11$ by inspection.

Finally, we check that the lowest rank of $Q_m(1,d_1,d_2, d_3, 0, \dots, 0)$ is at most $mn/2 - dm$ for $d = 1, 2$ and $m,n \geq 8d$.   The largest possible value of $d_3$, depending on $n \md{4}$, is given by:
\[ n = 0 \md{4} \implies (d_0, d_1, d_2,d_3) = (1,2,1, (n-8)/4). \]
\[ n = 1 \md{4} \implies (d_0, d_1, d_2,d_3) = (1,1,2, (n-9)/4). \]
\[ n = 2 \md{4} \implies (d_0, d_1, d_2,d_3) = (1,1,1, (n-6)/4). \]
\[ n = 3 \md{4} \implies (d_0, d_1, d_2,d_3) = (1,2,2, (n-11)/4). \]
Therefore, the largest possible value for the lowest rank $2d_1 + 6d_2 + 12 d_3$ is $3n-10$.  Now:
\[ mn/2 - dm \geq 4n - 8d \geq 3n + n - 8d \geq 3n \geq 3n - 10, \]
so $mn/2 - dm \geq 3n-10$ for $m,n \geq 8d$ as desired.
\end{proof}

\sse{Remark}  The lower bounds on $m$,$n$, and $r$ given above are certainly not the best possible.  One could improve the bounds on $m$ and $n$ by calculating the number of strictly unimodal products of five or more chains.  One could also improve the bound on $r$ by finding the overlaps among the strictly unimodal ranges for all products of three or more chains.   The lower bound on $r$ given in \cite{Za} is quadratic in $d$, while our lower bound is linear in $d$.  The best possible lower bound for $r$ is logarithmic in $d$, since it involves the inverse of the number of partitions of $r$.

 \end{document}